\documentclass[12pt,oneside,reqno]{amsart}
\usepackage{txfonts}
\usepackage{bbm}
\usepackage{amsmath}
\usepackage{graphicx}
\usepackage{mathrsfs}
\usepackage{stmaryrd}
\usepackage{amsfonts}
\usepackage{enumerate,amsmath,amssymb,amsthm}
\numberwithin{equation}{section}

\newcommand{\be}{\begin{eqnarray}}
\newcommand{\ee}{\end{eqnarray}}
\newcommand{\ce}{\begin{eqnarray*}}
\newcommand{\de}{\end{eqnarray*}}
\newtheorem{theorem}{Theorem}[section]
\newtheorem{lemma}[theorem]{Lemma}
\newtheorem{remark}[theorem]{Remark}
\newtheorem{definition}[theorem]{Definition}
\newtheorem{proposition}[theorem]{Proposition}
\newtheorem{Examples}[theorem]{Example}
\newtheorem{corollary}[theorem]{Corollary}

\def\eps{\varepsilon}

\def\e{\mathrm{e}}

\def\p{\partial}

\def\[{{\Big[}}
\def\]{{\Big]}}
\def\<{{\langle}}
\def\>{{\rangle}}
\def\({{\Big(}}
\def\){{\Big)}}

\def\bx{{\mathbf{x}}}

\def\dif{{\mathord{{\rm d}}}}

\def\no{\nonumber}
\def\={&\!\!=\!\!&}
\def\bt{\begin{theorem}}
\def\et{\end{theorem}}
\def\bl{\begin{lemma}}
\def\el{\end{lemma}}
\def\br{\begin{remark}}
\def\er{\end{remark}}
\def\bd{\begin{definition}}
\def\ed{\end{definition}}
\def\bp{\begin{proposition}}
\def\ep{\end{proposition}}
\def\bc{\begin{corollary}}
\def\ec{\end{corollary}}
\def\bx{\begin{Examples}}
\def\ex{\end{Examples}}

\def\cL{{\mathcal L}}

\def\mE{{\mathbb E}}

\def\mN{{\mathbb N}}

\def\mP{{\mathbb P}}
\def\mQ{{\mathbb Q}}
\def\mR{{\mathbb R}}
\def\mS{{\mathbb S}}

\def\sB{{\mathscr B}}

\def\sF{{\mathscr F}}

\def\sI{{\mathscr I}}

\def\sU{{\mathscr U}}

\def\geq{\geqslant}
\def\leq{\leqslant}

\allowdisplaybreaks

\begin{document}

\bf\title{Harnack inequalities for SDEs Driven by
Cylindrical\\ $\alpha$-stable processes}

\date{}
\author{Linlin Wang and Xicheng Zhang}

\address{Linlin Wang: School of Mathematics and Statistics, Wuhan University,
Wuhan, Hubei 430072, P.R.China\\
Email: Wanglinlin201296@gmail.com
}
 
\address{Xicheng Zhang: School of Mathematics and Statistics, Wuhan University,
Wuhan, Hubei 430072, P.R.China\\
Email: XichengZhang@gmail.com
 }

\dedicatory{}

\begin{abstract}
By using the coupling argument, we establish the Harnack and log-Harnack inequalites for stochastic differential equations with non-Lipschitz
drifts and driven by additive anisotropic subordinated Brownian motions (in particular, cylindrical $\alpha$-stable processes).
Moreover, the gradient estimate is also derived when the drift is Lipschitz continuous.

\end{abstract}

\maketitle

\rm

%\vspace{2mm}

\section{Introduction}

Consider the following stochastic differential equation (SDE):
\begin{align}
\dif X_t=b(X_t)\dif t+\dif L_t,\ \ X_0=x,\label{SDE}
\end{align}
where $L_t$ is a L\'evy process on some probability space $(\Omega,\sF,\mP)$ and $b:\mR^d\to\mR^d$ is continuous and satisfies
\begin{align}
(x_j-y_j)(b_j(x)-b_j(y))\leq |x_j-y_j|\rho(\|x-y\|_1),\ j=1,\cdots,d,\label{Cb}
\end{align}
for some $\rho\in\sU$. Here, $\|x\|_1:=\sum_j |x_j|$  and
\begin{align*}
\sU&:=\Big\{\mbox{$\rho:(0,\infty)\to(0,\infty)$ is a continuous non-decreasing}\\
&\qquad \mbox{and linear growth function with $\int_{0+}\frac{1}{\rho(s)}\dif s=+\infty$}\Big\}.
\end{align*}
Under the above assumptions, by Birhari's inequality (see Lemma \ref{Le21} below), it is easy to see that SDE (\ref{SDE}) has a unique solution denoted by $X_t(x)$, which defines
a Markov semigroup $(P_t)_{t\geq 0}$ by
$$
P_t f(x):=\mE f(X_t(x)),\ \ t\geq 0.
$$
where $f$ belongs to the class $\sB_b(\mR^d)$ of bounded measurable functions.

The aim of this work is to establish the Harnack inequalities for SDE (\ref{SDE}). The dimension-free Harnack inequality with power
was first introduced by F.Y. Wang in \cite{Wa1}, which can be used in the studies of heat kernel estimates, functional inequalities,
transportation-cost inequalities and properties of invariant measures, etc. (see \cite{Wa3}).
Up to now, the Harnack inequality with power and the log-Harnack inequality have been deeply studied
for stochastic (partial) differential equations driven by Brownian motions by using the coupling argument (cf. \cite{Wa3}).
In particular, we mention that by constructing a coupling with an unbounded time-dependent drift, F.Y. Wang in \cite{Wa2}
established the dimension-free Harnack inequalities for SDE (\ref{SDE}) driven by multiplicative Brownian noises 
under some monotonic conditions (see also \cite{Sh-Wa-Yu} and \cite{Zh2} for some extensions
with non-Lipschitz coefficients).
However, the corresponding results for SDEs driven by purely jump processes are very limited (see \cite{Sc-Wa, Sc-Sz-Wa} for the studies of L\'evy processes and
see \cite{O, Wa4, Wa5} for the studies of Ornstein-Uhlenbeck processes with jumps).
What is a coupling? Roughly speaking,
for a given stochastic curve starting from a fixed point (for example, the solution of an SDE),
we want to construct another stochastic curve starting from another fixed point (the solution of another SDE)
so that they can touch at a fixed time. In the case of jump diffusions,
since there are infinitely many jumps, it is usually hard to construct a coupling for the nonlinear equations.

In this work, we assume that $L_t$ takes the following form:
\begin{align}
L_t:=W_{S(t)}:=\Big(W^1_{S_1(t)},\cdots, W^d_{S_d(t)}\Big),\label{EE8}
\end{align}
where $W_t:=(W^1_t,\cdots, W^d_t)$ is a $d$-dimensional standard Brownian motion  on probability space $(\Omega,\sF,\mP)$,
and $S(t):=(S_1(t),\cdots, S_d(t))$ is an independent $d$-dimensional L\'evy process with each component $S_j(t)$ being a strictly positive subordinator with
Laplace transform given by
\begin{align}
\mE(\e^{-z\cdot S_t})=\exp\left\{-t\vartheta\cdot z+\int_{\mR^d_+}(\e^{-z\cdot u}-1)\nu_S(\dif u)\right\},\label{EW111}
\end{align}
where $\vartheta\in\mR^d_+$ and the L\'evy measure $\nu_S$ satisfies
$$
\int_{\mR^d_+}(1\wedge |u|)\nu_S(\dif u)<\infty.
$$
By easy calculations, one can see that the characteristic function of $L_t$ is given by
\begin{align}
\mE \mathrm{e}^{\mathrm{i} z\cdot L_t}=\exp\left\{-t\sum_k\vartheta_k|z_k|^2+t\int_{\mR^d}(\mathrm{e}^{\mathrm{i}z\cdot y}
-1-\mathrm{i}z\cdot y1_{|y|\leq 1})\nu_L(\dif y)\right\},\label{Ch}
\end{align}
where $\nu_L$ is the L\'evy measure given by
\begin{align}
\nu_L(\Gamma)=\int_{\mR^d_+}\left(\int_\Gamma\frac{(2\pi)^{-d/2}}{(u_1\cdots u_d)^{\frac{1}{2}}}
\mathrm{e}^{-(\frac{y_1^2}{2u_1}+\cdots+\frac{y_d^2}{2u_d})}\dif y_1\cdots\dif y_d\right)\nu_S(\dif u_1,\cdots,\dif u_d).\label{EW2}
\end{align}
Here we use the convention that if $u_i=0$ for some $i$, then the inner integral is calculated with respect to the degenerate Gaussian distribution.
In particular, $\nu_L$ may not be absolutely continuous with respect to the Lebesgue measure. However, obviously, $\nu_L$ is a symmetric measure.
Notice that for $\beta>0$,
$$
\int^\infty_0 u^{-1-\beta} \e^{-\frac{|y|^2}{2u}}\dif u=\Big(\tfrac{|y|^2}{2}\Big)^{-\beta}\Gamma(\beta),
$$
where $\Gamma$ is the usual Gamma function. If $\nu_S$ takes the following form:
$$
\nu_S(\dif u_1,\cdots,\dif u_d)=\int^\infty_0u^{-1-\frac{\alpha}{2}}\delta_u(\dif u_1)\cdots\delta_u(\dif u_d)\dif u,\ \  \alpha\in(0,2),
$$
then one sees that
$$
\nu_L(\dif y)=2^{-\frac{d+\alpha}{2}}(2\pi)^{-\frac{d}{2}}\Gamma(\tfrac{d+\alpha}{2})|y|^{-d-\alpha}\dif y_1\cdots\dif y_d.
$$
In this case, the generator of $L_t$ is given by
$$
\cL f(x):=\lim_{\eps\downarrow 0}\int_{|z|\geq\eps}[f(x+y)-f(x)]\nu_L(\dif y)=c_0(-\Delta)^{\frac{\alpha}{2}}f(x).
$$
If $\nu_S$ takes the following form:
$$
\nu_S(\dif u_1,\cdots,\dif u_d)=\sum_{i=1}^d\delta_0(\dif u_1)\cdots\delta_0(\dif u_{i-1})\Big(1_{u_i>0}u_i^{-1-\frac{\alpha_i}{2}}\Big)\dif u_i
\delta_0(\dif u_{i+1})\cdots\delta_0(\dif u_d), \ \  \alpha_i\in(0,2),
$$
then one sees that
$$
\nu_L(\dif y)=\sum_{i=1}^d2^{-\frac{1+\alpha_i}{2}}(2\pi)^{-\frac{1}{2}}\Gamma(\tfrac{1+\alpha_i}{2})\delta_0(\dif y_1)\cdots\delta_0(\dif y_{i-1})\Big(|y_i|^{-(1+\alpha_i)}\Big)\dif y_i
\delta_0(\dif y_{i+1})\cdots\delta_0(\dif y_d).
$$
In this case, we have
$$
\cL f(x)=\sum_{i=1}^dc_i(-\p_i^2)^{\frac{\alpha_i}{2}}f(x).
$$

Recently, X. Zhang in \cite{Zh1} used the time change and a smoothing approximation to derive a Bismut's type formula for SDE (\ref{SDE}) driven by subordinated Brownian motions
(see also \cite{Wa-Xu-Zh} for the extension of multiplicative noises by using a different approximation).
In \cite{Wa-Wa}, F.-Y. Wang and J. Wang used the same idea together with the coupling argument to derive the dimension-free Harnack inequalities for SDE (\ref{SDE}).
It should be emphasized that in these two works $L_t$ takes the same form:
$$
L_t=\Big(W^1_{S(t)},\cdots, W^d_{S(t)}\Big),
$$
where $S(t)$ is a one dimensional subordinator. In this case, $L_t$ is {\it isotropic} and the L\'evy measure of $L_t$ is {\it absolutely continuous} with respect to the Lebesgue measure.

In this paper, we are interested in the {\it anisotropic} case of (\ref{EE8}),  i.e,
each component may not jump simultaneously. In this case, the L\'evy measure of $L_t$ is {\it singular}.
We shall use the same arguments as in \cite{Zh1} and \cite{Wa-Wa} to prove the Harnack inequalities.
Since the L\'evy processes we are considering may have different scales
in different directions, the construction of the coupling will become more difficult,
and unfortunately, the price we have to pay is that the Harnack inequalities obtained below will be dimension-dependent.
Moreover, the time change is unapplicable since we have different clocks in different coordinates. Thus, it is so different from \cite{Wa-Wa} that we shall use
the Girsanov theorem for general continuous martingales rather than the one for Brownian motions.

Before stating our main result, we shall introduce the following functions: For $\rho\in\sU$,  define
\begin{align}
G_\rho(r):=\int^r_1\frac{1}{\rho(u)}\dif u,\ \ r>0,\label{ER88}
\end{align}
and for $T, r>0$,
\begin{align}
\Gamma_\rho(T,r):=r+T\rho\circ G^{-1}_\rho(G_\rho(r)+T).\label{ER8}
\end{align}
By definitions, these two functions are well-defined, and if $\rho(r)=C_0r$, then
\begin{align}
\Gamma_\rho(T,r):=(C_0 T\e^{C_0 T}+1) r.\label{EW9}
\end{align}
Our main result is:
\bt\label{Main}
Suppose (\ref{Cb}). We have the following conclusions:
\begin{enumerate}[(i)]
\item For any $T>0$ and strictly positive bounded measurable function $f$,
\begin{align}
P_T\log f(y)\leq\log P_T f(x)+\frac{1}{2}\Gamma^2_{\rho}(dT,\|x-y\|_1)\left(\sum_j\mE S_j(T)^{-1}\right),\quad x,y\in\mR^d.\label{ER5}
\end{align}
\item For any $T>0$, $p>1$ and bounded nonnegative measurable function $f$,
\begin{align}
(P_T f(y))^p\leq P_T f^p(x)\left(\mE\exp\left[\frac{p\Gamma^2_{\rho}(dT,\|x-y\|_1)}{2(p-1)^2}\sum_jS_j(T)^{-1}\right]\right)^{p-1}
,\quad x,y\in\mR^d.\label{ER55}
\end{align}
\end{enumerate}
\et

We shall prove this result in the next section. As a corollary, we have
\bt
Assume that $\|b\|_{\mathrm{Lip}}:=\sup_{x\not=y}\frac{\|b(x)-b(y)\|_1}{\|x-y\|_1}<\infty$. Then for any $T>0$ and bounded measurable function $f$,
\begin{align}
|\nabla P_Tf(x)|^2\leq [P_T f^2(x)-P_T f(x)^2](1+\|b\|_{\mathrm{Lip}}dT\e^{\|b\|_{\mathrm{Lip}}dT})^2d^2T^2\left(\sum_j\mE S_j(T)^{-1}\right),\label{ER4}
\end{align}
where
$$
|\nabla P_Tf(x)|:=\varlimsup_{y\to x}\frac{|P_Tf(x)-P_Tf(y)|}{\|x-y\|_1}.
$$
\et
\begin{proof}
It follows by (\ref{EW9}) with $\rho(r)=\|b\|_{\mathrm{Lip}}r$, (\ref{ER5}) and \cite[Proposition 2.3]{Ar-Th-Wa}.
\end{proof}
\br
Assume that $S_j(t)$ is an $\alpha_j/2$-subordinator, where $\alpha_j\in(0,2)$. It is well known that (see for example \cite[Proof of Theorem 1.1]{G}),
$$
\mE S_j(T)^{-1}\leq C_0T^{-\frac{2}{\alpha_j}},
$$
and if $\alpha_j\in(1,2)$, then for all $\lambda>0$ and $T>0$,
$$
\mE \e^{\lambda S_j(T)^{-1}}\leq 1+\left(\exp\left[\frac{C_1\lambda^{\frac{\alpha_j}{2(\alpha_j-1)}}}{T^{\frac{1}{\alpha_j-1}}}\right]-1\right)^{\frac{2(\alpha_j-1)}{\alpha_j}}
\leq\exp\left[\frac{C_2\lambda}{T^{\frac{2}{\alpha_j}}}+\frac{C_2\lambda^{\frac{\alpha_j}{2(\alpha_j-1)}}}{T^{\frac{1}{\alpha_j-1}}}\right].
$$
Applications of the above results are referred to \cite{Wa3}.
\er

\section{Proof of Main Theorem}

We need the following nonlinear Gronwall's inequality (cf. \cite[Lemma 2.1]{Zh-Zh}).
\bl\label{Le21}
(Bihari) Let $f:[0,\infty)\to[0,\infty)$ be a nonnegative measurable function. If for some $\rho\in\sU$,
\begin{align}
f(t)\leq f(0)+\int^t_0\rho(f(s))\dif s, ~~\forall t\geq 0,\label{e8}
\end{align}
then
\begin{align}
f(t)\leq G_\rho^{-1}\left(G_\rho(f(0))+t\right),\ \ \forall t\geq 0, \label{e12}
\end{align}
where $G_\rho$ is defined by (\ref{ER88}). In particular, if $f(0)=0$, then $f(t)\equiv0$.
\el
\begin{proof}
Set
$$
h(t):=f(0)+\int^t_0\rho(f(s))\dif s.
$$
Then $f(t)\leq h(t)$, and by the non-decrease of $\rho$, we have
$$
\rho(f(t))\leq \rho(h(t)).
$$
By the usual differential formula,  we have
\begin{align*}
G_\rho(h(t))=G_\rho(f(0))+\int^t_0G_\rho'(h(s))\dif h(s)
= G_\rho(f(0))+\int^t_0\frac{\rho(f(s))}{\rho(h(s))}\dif s
\leq G_\rho(f(0))+t,
\end{align*}
which in turn implies (\ref{e12}). If $f(0)=0$, by $G_\rho(0)=-\infty$ and $G^{-1}_\rho(-\infty)=0$, we get $f(t)\equiv 0$.
\end{proof}

Let $\mS^d$ be the space of all right continuous functions from $[0,\infty)\to[0,\infty)^d$ and having left hand limits, which satisfies that for each $\ell=(\ell_1,\cdots,\ell_d)\in\mS^d$,
$s\mapsto\ell_j(s)$ is increasing with $\ell_j(0)=0$.
Let $X^\ell_t(x)$ solve the following SDE:
\begin{align}
\dif X^\ell_t=b(X^\ell_t)\dif t+\dif W_{\ell(t)},\ \ X^\ell_0=x.
\end{align}
We first prove that
\bt\label{Th2}
Suppose (\ref{Cb}).  We have the following conclusions:
\begin{enumerate}[(i)]
\item For any $T>0$ and strictly positive bounded measurable function $f$,
\begin{align}
\mE\log f(X^\ell_T(y))\leq\log \mE f(X^\ell_T(x))+\frac{1}{2}\Gamma^2_{\rho}(dT,\|x-y\|_1)\left(\sum_j\ell_j(T)^{-1}\right),\quad x,y\in\mR^d.\label{ER9}
\end{align}
\item For any $T>0$, $p>1$ and bounded nonnegative measurable function $f$,
\begin{align}
(\mE f(X^\ell_T(y)))^p\leq \mE(f(X^\ell_T(x)))^p\left(\exp\left[\frac{p\Gamma^2_{\rho}(dT,\|x-y\|_1)}{2(p-1)^2}\sum_j\ell_j(T)^{-1}
\right]\right)^{p-1},\quad x,y\in\mR^d.\label{ER10}
\end{align}
\end{enumerate}
\et
\subsection{$\ell$ being absolutely continuous and strictly increasing}
Let $\ell=(\ell_1,\cdots,\ell_d)$ with $\ell_j$ being an absolutely continuous and strictly increasing function from $[0,\infty)$ to $[0,\infty)$ with $\ell_j(0)=0$.
Fix $T>0$ and $x=(x_1,\cdots,x_d), y=(y_1,\cdots,y_d)\in\mR^d$. Consider the following coupled SDE:
\begin{align}
\left\{\label{Co}
\begin{aligned}
X^j_t&=x_j+\int^t_0b_j(X_s)\dif s+W^j_{\ell_j(t)},\quad \ j=1,\cdots, d,\\
Y^j_t&=y_j+\int^t_0b_j(Y_s)\dif s+W^j_{\ell_j(t)}+\kappa_T\int^{t\wedge\tau_j}_0\frac{X^j_s-Y^j_s}{|X^j_s-Y^j_s|}\dif \ell_j(s),
\end{aligned}
\right.
\end{align}
where $\tau_j:=\inf\{t\geq 0: X^j_t=Y^j_t\}$ and $\kappa_T$ will be chosen below so that
$$
X_T=Y_T.
$$
Below, we shall use the following filtration: for $\ell\in\mS^d$,
$$
\sF^\ell_t:=\sigma\Big\{W^j_s: s\leq {\ell_j(t)},\ j=1,\cdots,d\Big\}.
$$
\bl\label{Le24}
Under (\ref{Cb}), there exists a unique pair of continuous ($\sF^\ell_t$)-adapted process $(X,Y)$ solving equation (\ref{Co}) such that for each $j=1,\cdots,d$,
\begin{align}
X^j_t\not=Y^j_t,\ \ t<\tau_j\mbox{ and } \ X^j_t=Y^j_t, \ \ t\geq\tau_j.\label{St}
\end{align}
\el
\begin{proof}
Roughly speaking, condition (\ref{St}) means that the components $X^j_t$ and $Y^j_t$ will go together after they meet at time $\tau_j$.
To construct the solution, we consider the following SDE:
$$
\left\{
\begin{aligned}
X^j_t&=x_j+\int^t_0b_j(X_s)\dif s+W^j_{\ell_j(t)},\quad \ j=1,\cdots, d,\\
\widetilde Y^j_t&=y_j+\int^t_0b_j(\widetilde Y_s)\dif s+W^j_{\ell_j(t)}
+\kappa_T\int^t_01_{\{X^j_s\not=\widetilde Y^j_s\}}\frac{X^j_s-\widetilde Y^j_s}{|X^j_s-\widetilde Y^j_s|}\dif \ell_j(s).
\end{aligned}
\right.
$$
Since $1_{\{x_j\not= y_j\}}\frac{x_j-y_j}{|x_j-y_j|}, j=1,\cdots,d$  are locally Lipschitz on
$\{(x,y): x_j\not= y_j,\ \forall j=1,\cdots,d\}$, this equation can be uniquely solved up to the time $\widetilde\tau_1$:
$$
\widetilde\tau_1:=\inf\Big\{t\geq 0: X^j_t=\widetilde Y^j_t, \exists j=1,\cdots,d\Big\}.
$$
Below we use the convention:
$$
X^j_\infty(\omega)=\widetilde Y^j_\infty(\omega)=\infty.
$$
We define
$$
N_1(\omega):=\Big\{j\in\{1,\cdots,d\}: X^j_{\widetilde\tau_1(\omega)}(\omega)=\widetilde Y^j_{\widetilde\tau_1(\omega)}(\omega)\Big\}
$$
and
$$
Y^j_t(\omega):=
\left\{
\begin{aligned}
&\widetilde Y^j_t(\omega),\ \ t<\widetilde\tau_1(\omega),\ j\in \{1,\cdots,d\};\\
&X^j_t(\omega),\ \ t\geq\widetilde\tau_1(\omega),\ \ j\in N_1(\omega).
\end{aligned}
\right.
$$
Next we consider the following SDE:
$$
\left\{
\begin{aligned}
X^j_t&=X^j_{\widetilde\tau_1}+\int^t_{\widetilde\tau_1}b_j(X_s)\dif s+W^j_{\ell_j(t)}-W^j_{\ell_j(\widetilde\tau_1)},\ \ j\notin N_1(\omega),\\
\widetilde Y^j_t&=Y^j_{\widetilde\tau_1}+\int^t_{\widetilde\tau_1}b_j((X^i_s)_{i\in N_1}, (\widetilde Y^i_s)_{i\notin N_1})\dif s+W^j_{\ell_j(t)}-W^j_{\ell_j(\widetilde\tau_1)}\\
&\qquad+\kappa_T\int^t_{\widetilde\tau_1}1_{\{X^j_s\not=\widetilde Y^j_s\}}\frac{X^j_s-\widetilde Y^j_s}{|X^j_s-\widetilde Y^j_s|}\dif \ell_j(s).
\end{aligned}
\right.
$$
This equation can be uniquely solved up to the time
$$
\widetilde\tau_2:=\inf\Big\{t\geq \widetilde\tau_1: X^j_t=\widetilde Y^j_t, \exists j\notin N_1\Big\}.
$$
We define
$$
N_2(\omega):=\Big\{j\notin N_1(\omega): X^j_{\widetilde\tau_2(\omega)}(\omega)=\widetilde Y^j_{\widetilde\tau_2(\omega)}(\omega)\Big\}
$$
and
$$
Y^j_t(\omega):=
\left\{
\begin{aligned}
&\widetilde Y^j_t(\omega),\ \ t<\widetilde\tau_2(\omega),\ j\notin N_1(\omega);\\
&X^j_t(\omega),\ \ t\geq\widetilde\tau_2(\omega),\ \ j\in N_2(\omega).
\end{aligned}
\right.
$$
Proceeding this construction at most $d$-times, we obtain a unique pair of $(X,Y)$ solving equation (\ref{Co})  and satisfying (\ref{St}).
Let $\sI_d$ be the set of all subsets of $\{1,\cdots,d\}$.
Then $N_1,N_2,\cdots$ can be regarded as random variables in $\sI_d$. Since $\sI_d$ has only finitely many elements and the above construction has at most $d$-steps,
if we restrict our consideration on $[0,t]$,  it is easy to see that $Y_t$ is ($\sF^\ell_t$)-measurable.
\end{proof}

Now, set
$$
Z^j_t:=X^j_t-Y^j_t,\ \ j=1,\cdots,d.
$$
By (\ref{Co}) we have
$$
Z^j_t=Z^j_0+\int^t_0(b_j(X_s)-b_j(Y_s))\dif s-\kappa_T\int^{t\wedge\tau_j}_0\frac{Z^j_s}{|Z^j_s|}\dif\ell_j(s).
$$
By the differential formula, we have
\begin{align*}
|Z^j_t|&=|Z^j_0|+\int^t_0|Z^j_s|^{-1}Z^j_s(b_j(X_s)-b_j(Y_s))\dif s-\kappa_T\ell_j(t\wedge\tau_j).
\end{align*}
By (\ref{Cb}), we obtain
$$
\|Z_t\|_1+\kappa_T\sum_j\ell_j(t\wedge\tau_j)\leq \|Z_0\|_1+d\int^t_0\rho(\|Z_s\|_1)\dif s.
$$
By Lemma \ref{Le21}, we obtain
$$
\|Z_t\|_1\leq G^{-1}_\rho(G_\rho(\|Z_0\|_1)+dt),\ \ t\geq 0,
$$
where $G_\rho(r)$ is defined by (\ref{ER88}) and
\begin{align}
\kappa_T\sum_{j=1}^d\ell_j(t\wedge\tau_j)\leq \|Z_0\|_1+dt\rho\circ G^{-1}_\rho(G_\rho(\|Z_0\|_1)+dt)=\Gamma_\rho(dt,\|Z_0\|_1).\label{EE1}
\end{align}
For any $\eps\in(0,1)$, if we choose
\begin{align}
\kappa_T:=\frac{\Gamma_\rho(dT,\|x-y\|_1)}{\ell_1(\eps T)\wedge\cdots\wedge\ell_d(\eps T)},\label{EE4}
\end{align}
then by (\ref{EE1}),
\begin{align}
\sum_{j}\ell_j(T\wedge\tau_j)\leq\ell_1(\eps T)\wedge\cdots\wedge\ell_d(\eps T).\label{EE2}
\end{align}
In particular,
$$
\ell_j(T\wedge\tau_j)\leq\ell_j(\eps T),\ \ \forall j=1,\cdots, d,
$$
which, by the strict increase of $\ell_j$ and $\eps<1$, implies that
\begin{align}
\tau_j\leq T,\ \ \forall j=1,\cdots, d.\label{EE3}
\end{align}
Hence, by (\ref{St}) one has
\begin{align}
X_T=Y_T.\label{ER7}
\end{align}

Now define
$$
H^j_s:=\kappa_T1_{\{s<\tau_j\}}\frac{X^j_s-Y^j_s}{|X^j_s-Y^j_s|}
$$
and
$$
M_t:=-\sum_j\int^t_0H^j_s\dif W^j_{\ell_j(s)}.
$$
Clearly, $M_t$ is a continuous ($\sF^\ell_t$)-martingale with
\begin{align}
\<M\>_t=\sum_j\int^t_0(H^j_s)^2\dif\ell_j(s)=\kappa_T^2\sum_j\ell_j(t\wedge\tau_j).\label{EE5}
\end{align}
Notice that by (\ref{EE4}), (\ref{EE2}) and (\ref{EE3}),
\begin{align}
\<M\>_\infty=\kappa_T^2\sum_j\ell_j(\tau_j)=\kappa_T^2\sum_j\ell_j(T\wedge\tau_j)\leq\frac{\Gamma^2_\rho(dT,\|x-y\|_1)}{\ell_1(\eps T)\wedge\cdots\wedge\ell_d(\eps T)}.\label{ER6}
\end{align}
Thus,  if we let
$$
R:=\exp\left\{M_\infty-\tfrac{1}{2}\<M\>_\infty\right\},
$$
then by Novikov's criterion,
$$
\mE R=1.
$$
Set
$$
\dif\mQ:=R\dif\mP,\ \ \widetilde W^j_{\ell_j(t)}:=W^j_{\ell_j(t)}+\int^t_0H^j_s\dif \ell_j(s).
$$
By Girsanov's theorem (cf. \cite{Re-Yo}), $\widetilde W_{\ell(t)}:=\big(\widetilde W^1_{\ell_1(t)}, \cdots,\widetilde W^d_{\ell_d(t)}\big)$
is a $d$-dimensional ($\sF^\ell_t$)-martingale under $\mQ$ and
\begin{align*}
\<\widetilde W^i_{\ell_i}, \widetilde W^j_{\ell_j}\>^\mQ_t=\<W^i_{\ell_i}, W^j_{\ell_j}\>^\mP_t=1_{i=j}\ell_i(t).
\end{align*}
For any $\eta=(\eta_1,\cdots,\eta_d)$, in view of
$$
\<\eta\cdot \widetilde W_\ell\>^\mQ_t=\sum_j\eta_j^2\ell_j(t),
$$
we have
$$
\exp\left\{\mathrm{i}\eta\cdot \widetilde W_{\ell(t)}+\frac{1}{2}\sum_j\eta_j^2\ell_j(t)\right\}\mbox{ is a complex ($\sF^\ell_t$)-martingale under $\mQ$.}
$$
Hence,
$$
\mE^\mQ\left[\exp\left\{\mathrm{i}\eta\cdot (\widetilde W_{\ell(t)}-\widetilde W_{\ell(s)})\right\}|\sF^\ell_s\right]=\exp\left\{\frac{1}{2}\sum_j\eta_j^2(\ell_j(t)-\ell_j(s))\right\},
$$
which implies that the law of $\widetilde W_\ell$ under $\mQ$ is the same as that of $W_{\ell}$ under $\mP$. Since
$$
Y^j_t=y_j+\int^t_0 b_j(Y_s)\dif s+\widetilde W^j_{\ell_j(t)},\ \ j=1,\cdots,d,
$$
we also have
\begin{align}
\mbox{the law of $X_t(y)$ under $\mP$ is the same as that of $Y_t(y)$ under $\mQ$.}\label{La}
\end{align}

\begin{proof}[Proof of Theorem \ref{Th2}:] Below, we shall drop the superscript $\ell$ in $X^\ell$ for simplicity of notations.

(i) Let $\mu$ be a probability measure on $\mR^d$. Notice the following fact that for any $g_1,g_2\geq 0$ with $\mu(g_1)=1$,
$$
\mu(g_1g_2)\leq \log\mu(\e^{g_2})+\mu(g_1\log g_1),
$$
where $\mu(\cdot)$ denotes the expectation with respect to $\mu$. By (\ref{La}) and (\ref{ER7}), we have
\begin{align}
\mE\log f(X_T(y))&=\mE[R\log f(Y_T(y))]=\mE[R\log f(X_T(x))]\no\\
&\leq \log\mE f(X_T(x))+\mE[R\log R].\label{ER1}
\end{align}
By the definition of $R$, we have
\begin{align*}
\mE (R\log R)&=\mE^\mQ\left\{M_\infty-\tfrac{1}{2}\<M\>_\infty\right\}\\
&=\sum_j\mE^\mQ\left\{-\int^\infty_0H^j_s\dif \widetilde W^j_{\ell_j(s)}+\frac{1}{2}\int^\infty_0(H^j_s)^2\dif\ell_j(s)\right\}\\
&=\frac{1}{2}\mE^\mQ\<M\>_\infty\stackrel{(\ref{ER6})}{\leq}\frac{\Gamma^2_\rho(dT,\|x-y\|_1)}{2\ell_1(\eps T)\wedge\cdots\wedge\ell_d(\eps T)}\\
&\leq\frac{1}{2}\Gamma^2_\rho(dT,\|x-y\|_1)\sum_j\ell_j(\eps T)^{-1},
\end{align*}
which, together with (\ref{ER1}) and letting $\eps\uparrow 1$, gives (\ref{ER9}).

(ii) For $p>1$, let $q:=\frac{p}{p-1}$. By (\ref{La}), (\ref{ER7}) again and H\"older's inequality, we have
\begin{align}
(\mE f(X_T(y)))^p=(\mE (R f(X_T(x))))^p\leq \left(\mE R^q\right)^{p-1}\mE(f(X_T(x)))^p.\label{ER2}
\end{align}
On the other hand, by the definition of $R$ and Novikov's criterion, we also have
\begin{align*}
\mE R^q&=\mE \exp\left\{qM_\infty-\tfrac{q}{2}\<M\>_\infty\right\}\\
&=\mE \exp\left\{qM_\infty-\tfrac{1}{2}\<qM\>_\infty+\tfrac{q^2-q}{2}\<M\>_\infty\right\}\\
&\stackrel{(\ref{ER6})}{\leq}\exp\left\{\frac{(q^2-q)\Gamma^2_\rho(dT,\|x-y\|_1)}{2\ell_1(\eps T)\wedge\cdots\wedge\ell_d(\eps T)}\right\},
\end{align*}
which, together with (\ref{ER2})  and letting $\eps\uparrow1$, implies (\ref{ER10}).
\end{proof}
\subsection{General $\ell$}
For $n\in\mN$, we define
$$
\ell^n_j(t):=n\int^{t+\frac{1}{n}}_t\ell_j(s)\dif s+\frac{t}{n}.
$$
Clearly, $t\mapsto\ell^n_j(t)$ is absolutely continuous and strictly increasing and
\begin{align}
\ell^n_j(t)\downarrow \ell_j(t),\ \ n\to\infty.\label{KJ4}
\end{align}
Let $X^{\ell^n}_t$ solve the following SDE:
$$
X^{\ell^n}_t=x+\int^t_0b(X^{\ell^n}_s)\dif s+W_{\ell^n(t)}.
$$
\bl\label{Le43}
Assume that $b$ is locally Lipschitz and linear growth, then for each $t\geq 0$, $X^{\ell^n}_t$ converges to $X^\ell_t$ in probability.
\el
\begin{proof}
Since $b$ is linear growth, it is easy to prove that for any $T>0$,
\begin{align}
\sup_{n\in\mN}\mE\left(\sup_{t\in[0,T]}\|X^{\ell^n}_t\|_1\right)<\infty.\label{KJ5}
\end{align}
For $R>0$, define
$$
\tau^n_R:=\inf\{t\geq 0: \|X^{\ell^n}_t\|_1\vee\|X^{\ell}_t\|_1>R\}.
$$
For any $t<\tau^n_R$, we have
\begin{align*}
\|X^{\ell^n}_t-X^\ell_t\|_1&\leq\int^t_0\|b(X^{\ell^n}_s)-b(X^{\ell}_s)\|_1\dif s+\|W_{\ell^n(t)}-W_{\ell(t)}\|_1\\
&\leq C_R\int^t_0\|X^{\ell^n}_s-X^{\ell}_s\|_1\dif s+\|W_{\ell^n(t)}-W_{\ell(t)}\|_1,
\end{align*}
which yields by Gronwall's inequality that for $t<\tau^n_R$,
$$
\|X^{\ell^n}_t-X^\ell_t\|_1\leq \|W_{\ell^n(t)}-W_{\ell(t)}\|_1+C_R\e^{t C_R}\int^t_0\|W_{\ell^n(s)}-W_{\ell(s)}\|_1\dif s.
$$
Now, for any $\eps>0$, by Chebyshev's inequality we have
\begin{align*}
\mP(\|X^{\ell^n}_t-X^\ell_t\|_1\geq\eps)&\leq\mP(t\geq \tau^n_R)+\mP(\|X^{\ell^n}_t-X^\ell_t\|_1\geq\eps; t<\tau^n_R)\\
&\leq\mP\left\{\sup_{s\in[0,t]}\|X^{\ell^n}_s\|_1\vee\|X^{\ell}_s\|_1>R\right\}+\mP\Big(\|W_{\ell^n(t)}-W_{\ell(t)}\|_1\geq\tfrac{\eps}{2}\Big)\\
&\quad+\mP\left\{C_R\e^{t C_R}\int^t_0\|W_{\ell^n(s)}-W_{\ell(s)}\|_1\dif s\geq\frac{\eps}{2}\right\}\\
&\leq\frac{1}{R}\mE\left(\sup_{s\in[0,t]}\|X^{\ell^n}_s\|_1\vee\|X^{\ell}_s\|_1\right)+\frac{2}{\eps}\mE\Big(\|W_{\ell^n(t)}-W_{\ell(t)}\|_1\Big)\\
&\quad+\frac{2}{\eps}\mE\left(C_R\e^{t C_R}\int^t_0\|W_{\ell^n(s)}-W_{\ell(s)}\|_1\dif s\right).
\end{align*}
First letting $n\to\infty$ and then $R\to\infty$, by (\ref{KJ4}) and (\ref{KJ5}), we obtain the desired convergence.
\end{proof}

Let $\rho$ be a nonnegative $C^\infty$-function on $\mR^d$ with support in the unit ball and
$$
\int_{\mR^d}\rho(x)\dif x=1.
$$
For $\eps\in(0,1)$, write
\begin{align*}
\rho_\eps(x):=\eps^{-d}\rho(x/\eps)
\end{align*}
and
\begin{align*}
b_\eps(x):=\int_{{\mR}^d}b(z)\rho_\eps(x-z)\dif z=\int_{{\mR}^d}b(x-z)\rho_\eps(z)\dif z.
\end{align*}
It is easy to see that for each $x\in\mR^d$,
\begin{align}
\lim_{\eps\downarrow 0}b_\eps(x)=b(x),\label{KJ2}
\end{align}
and by (\ref{Cb}),
\begin{align}\label{Cbd}
(x_j-y_j)(b_\eps^j(x)-b_\eps^j(y))\leq|x_j-y_j|\rho(\|x-y\|_1).
\end{align}
Let $X^\eps_t$ solve the following equation:
\begin{align}\label{SDE2}
X^\eps_t=x_0+\int_0^tb_\eps(X^\eps_s)\dif s+W_{\ell(t)}.
\end{align}
\bl\label{Le44}
For each $\omega\in\Omega$ and $t\geq 0$, we have
$$
\lim_{\eps\downarrow 0}|X^\eps_t(\omega)-X^\ell_t(\omega)|=0.
$$
\el
\begin{proof}
Since $b_\eps$ is uniformly linear growth with respect to $\eps$, it is easy to prove that for each $\omega$ and $T>0$,
\begin{align}
\sup_{\eps\in(0,1)}\sup_{t\in[0,T]}\|X_{t}^{\eps}(\omega)\|_1<\infty.\label{KJ1}
\end{align}
Denote
\begin{align*}
Z_t^{\eps,j}:=X_t^{\eps,j}-X_t^{\ell,j},\quad j=1,\cdots,d,
\end{align*}
then
\begin{align*}
Z_t^{\eps,j}=\int_0^t(b_{\eps}^j(X_s^{\eps})-b^j(X_s^{\ell}))\dif s.
\end{align*}
By the differential formula and (\ref{Cbd}), we obtain
\begin{align*}
|Z_t^{\eps,j}|
&=\int_0^t|Z_s^{\eps,j}|^{-1}Z_s^{\eps,j}(b_{\eps}^j(X_{s}^{\eps})-b^j(X_s^{\ell}))\dif s
\\&\leq\int_0^t\rho(\|Z_s^\eps\|_1)\dif s
+\int_0^t|b_{\eps}^j(X_s^{\ell})-b^j(X_s^{\ell})|\dif s.
\end{align*}
Thus,
\begin{align*}
\|Z_t^{\eps}\|_1\leq d\int_0^t\rho(\|Z_s^\eps\|_1)\dif s+\int_0^t\|b_{\eps}(X_s^{\ell})-b(X_s^{\ell})\|_1\dif s.
\end{align*}
Taking limits for both sides and by Fatou's lemma and (\ref{KJ1}), (\ref{KJ2}), we arrive at
\begin{align*}
\varlimsup_{\eps\downarrow 0}\|Z_t^\eps\|_1\leq d\int_0^t\rho(\varlimsup_{\delta\downarrow 0}\|Z_s^\eps\|_1)\dif s,
\end{align*}
which yields the desired result by Lemma \ref{Le21}.
\end{proof}

\begin{proof}[Proof of Theorem \ref{Th2}:] Since we have proven in the previous section that
\begin{align*}
\mE\log f(X^{\ell^n}_T(y))\leq\log \mE f(X^{\ell^n}_T(x))+\frac{1}{2}\Gamma^2_\rho(dT,\|x-y\|_1)\left(\sum_j\ell^n_j(T)^{-1}\right),
\end{align*}
and
\begin{align*}
|\mE f(X^{\ell^n}_T(y))|^p\leq \mE|f(X^{\ell^n}_T(x))|^p\left(\exp\left[\frac{p\Gamma^2_\rho(dT,\|x-y\|_1)}{2(p-1)^2}\sum_j\ell^n_j(T)^{-1}\right]\right)^{p-1},
\end{align*}
by Lemma \ref{Le43} and taking limit $n\to\infty$, we obtain the desired estimates for the general $\ell$ when $b$ is locally Lipschitz and linear growth.
Now that $b_{\eps}(x)$ is locally Lipschitz and linear growth,
the claimed inequalities hold for SDE (\ref{SDE2}). Hence, by Lemma \ref{Le44} and letting $\eps\downarrow0$, we get the desired estimates for the general $\ell$
when $b$ satisfies (\ref{Cb}).
\end{proof}

\vspace{5mm}

\begin{proof}[Proof of Theorem \ref{Main}:]
Let $X_t(x)$ be the solution of SDE (\ref{SDE}). Clearly, we have
$$
P_Tf(x)=\mE f(X_T(x))=\mE(\mE f(X^\ell_t(x))|_{\ell=S}).
$$
Basing on this, (\ref{ER5}) follows by (\ref{ER9}) and Jensen's inequality, and (\ref{ER55}) follows by (\ref{ER10}) and H\"older's inequality.
\end{proof}

{\bf Acknowledgements:}

The authors are grateful to Professors Feng-Yu Wang and Jian Wang for their quite useful suggestions.
This work is supported by NSFs of China (No. 11271294).

\end{document}